\documentclass{amsart}
\usepackage{amssymb,amsmath,amsfonts}

\newcommand{\C}{\mathbb C}

\newcommand{\N}{\mathbb N}
\newcommand{\Z}{\mathbb Z}
\newcommand{\T}{\otimes}
\newcommand{\Cn}{\C[e_0,\ldots,e_{n-1}]}
\newcommand{\an}{a_1,\ldots,a_n}

\newcommand{\ld}{\ldots}
\newcommand{\hk}{\hookrightarrow}

\newcommand{\Sum}{\sum\limits}
\newcommand{\Prod}{\prod\limits}

\newcommand{\nc}{\newcommand}
\nc{\g}{\mathfrak g}
\nc{\h}{\mathfrak h}
\nc{\n}{\mathfrak n_+}

\nc{\pin}{\pi_1,\ld,\pi_n}
\nc{\zn}{z_1,\ld,z_n}
\nc{\slt}{\mathfrak{sl}_2}
\nc{\slth}{\widehat{\slt}}
\nc{\Tpi}{\pi_1\T\ld\T\pi_n}
\nc{\fpi}{\pi_1*\ld*\pi_n}
\nc{\A}{\mathfrak A}
\nc{\yn}{y_1,\ld,y_n}
\nc{\pz}{\phi_Z}
\nc{\al}{\alpha}
\nc{\be}{\beta}
\nc{\CT}{(\C^2)^{\T n}}
\nc{\Ca}{\C^{a_1}\T\ld\T\C^{a_n}}
\nc{\veps}{\varepsilon}
\nc{\ik}{i_1,\ld,i_k}
\nc{\ws}{\widetilde\sigma}
\nc{\W}{\widetilde W}

\nc{\bs}{b_1,\ldots,b_s}

\nc{\qb}[2]{{\genfrac{[}{]}{0pt}{0}{#1}{#2}}_q}
\nc{\gf}[2]{\genfrac{}{}{0pt}{}{#1}{#2}}
\nc{\qf}[1]{({#1})_q!}

\nc{\codim}{{\mathop{\rm codim}}}
\nc{\Id}{{\mathop{\rm Id}}}
\nc{\Tr}{{\mathop{\rm Tr}}}
\nc{\ch}{{\mathop{\rm ch}}}
\nc{\V}{{\mathop{\rm V}}}
\nc{\const}{{\mathop{\rm const}}}
\nc{\sltc}{\slt\otimes\C[t]}

\newtheorem{rem}{Remark}[section]
\newtheorem{opr}{Definition}[section]
\newtheorem{theorem}{Theorem}[section]
\newtheorem{utv}{Statement}[section]
\newtheorem{lem}{Lemma}[section]
\newtheorem{cor}{Corollary}[section]
\newtheorem{prop}{Proposition}[section]

\begin{document}
\pagestyle{plain}
\markboth{B.Feigin and E.Feigin}
{Integrable $\slth$-modules as infinite tensor products}
\title[Integrable $\slth$-modules as infinite tensor products]
{Integrable $\slth$-modules as infinite tensor products}
\author{B.Feigin and E.Feigin}
\address{B.F.: Landau Institute for Theoretical Physics, Russia,
Chernogolovka, 142432}
\email {feigin@mccme.ru }
\address{E.F.:Independent University of Moscow, Russia, Moscow,
Bol'shoi Vlas'evskii per.,7.}
\email {evgfeig@mccme.ru }
\date{}
\begin{abstract}
In this paper we develop an ideas and methods from
\cite{mi}. Using the fusion product of the representations of the
Lie algebra $\slt$,
we construct a set of the integrable highest weight $\slth$-modules $L^D$,
depending on the vector $D\in \N^{k+1}$. In a special cases of $D$
our modules are isomorphic to the irreducible $\slth$-modules $L_{i,k}$.
We construct a basis of the $L^D$ and study the decomposition of
$L^D$ on the irreducible components. We also write down a formulas
for the characters of $L^D$.
\end{abstract}

\maketitle
\section{Introduction}
In this paper we introduce a way to construct $\slth$-modules
as a limit of the finite-dimensional representations of the polynomial
current Lie algebra. We use the notion of the fusion product of the
finite-dimensional representations of $\slt$, introduced in
\cite{fusion} and studied in \cite{mi}. Let us recall main definitions.

Let $\slt\otimes\C[t]$ be the polynomial current algebra, $e,f,h$ -- a
standard $\slt$-basis.
Let $\pi$ be an $\slt$-module, $z\in\C$. Define an $\slt\otimes\C[t]$-action
on the $\pi$: $(x\otimes t^i)\cdot v=z^i x\cdot v,\ x\in\slt,v\in\pi$.
Denote this $\slt\otimes\C[t]$-module by $\pi(z)$.

Let $\pin$ be an irreducible representations of $\slt$, $\dim\pi_i=a_i$.
Let
$Z=(z_1,\ld,z_n)\in\C^n$, all $z_i$ are pairwise distinct.
Consider the tensor product of the $\slt\otimes\C[t]$-modules
$\pi_1(z_1)\otimes\ld\otimes\pi_n(z_n)$.
If $w$ is a product of the highest vectors of $\pi_i$, then
$\pi_1(z_1)\otimes\ld\otimes\pi_n(z_n)$ is a cyclic
$\slt\otimes\C[t]$-module with a cyclic vector $w$,
$$\pi_1(z_1)\otimes\ld\otimes\pi_n(z_n)=U(\slt\otimes\C[t])w.$$
Denote by $\fpi$ an adjoint graded module with respect to the filtration
by the $t$-degree on $U(\slt\otimes\C[t])$. We call this module the
fusion product of $\pin$.

In the paper we prove the existence of injective homomorphisms
between different fusion products. Let $a_1\le\ld\le a_n$.
Then, if $a_i=a_{i+1}$, we construct an injection
\begin{equation}
\label{inj1}
\C^{a_1}*\ld *\C^{a_{i-1}}*\C^{a_{i+2}}*\ld *\C^{a_n}\hk
\C^{a_1}*\ld *\C^{a_n}.
\end{equation}
Moreover, there exists an isomorphism of the $\sltc$-modules
\begin{multline}
\label{quot1}
(\C^{a_1}*\ld *\C^{a_n})/
(\C^{a_1}*\ld *\C^{a_{i-1}}*\C^{a_{i+2}}*\ld *\C^{a_n})\cong\\
\cong
\C^{a_1}*\ld *\C^{a_{i-1}}*\C^{a_i-1}*
\C^{a_{i+1}+1}*\C^{a_{i+2}}*\ld *\C^{a_n}.
\end{multline}
(Note that the similar formulas were used by Shilling and Warnaar to
construct the supernomials, see \cite{sup,mult,skost}).
Now, suppose $a_1<a_2$. Then we construct an injection
\begin{equation}
\label{inj2}
\C^{a_2-a_1+1}*\C^{a_3}*\ld *\C^{a_n}\hk
\C^{a_1}*\ld *\C^{a_n}
\end{equation}
with a quotient:
\begin{multline}
\label{quot2}
(\C^{a_1}*\ld *\C^{a_n})/
(\C^{a_2-a_1+1}*\C^{a_3}*\ld *\C^{a_n})\cong \\ \cong
\C^{a_1-1}*\C^{a_2+1}*\C^{a_3}*\ld *\C^{a_n}.
\end{multline}

In fact, the more general statement is true. Let $i<j$. Then there
exist submodules $S_{i,j}(A)=S_{i,j}(a_1,\ld,a_n)$ of the fusion product
of $\C^{a_i}$, such that the following holds:
\begin{multline}
\C^{a_1}*\ld *\C^{a_n}/S_{i,j}(A)\cong\\ \cong
\C^{a_1}*\ld *\C^{a_{i-1}}*\C^{a_i-1}*\C^{a_{i+1}}*\ld
*\C^{a_{j-1}}*\C^{a_j+1}*\C^{a_{j+1}}*\ld *\C^{a_n}.
\end{multline}
Looking at the formula (\ref{quot2}), one can expect $S_{i,j}(A)$ to be equal
to the following fusion product:
\begin{equation}
S_{i,j}(A)\cong
\C^{a_j-a_i+1}*\C^{a_1}*\ld*\C^{a_{i-1}}*\C^{a_{i+1}}*\ld *
\C^{a_{j-1}}*\C^{a_{j+1}}*\ld *\C^{a_n}.
\end{equation}
But that is not true in the case of the general $i,j$. However, the
structure of the $\sltc$-modules $S_{i,j}(A)$ is very interesting.
We do not study it in this paper, but we hope to return to this topic later.

One can consider the formulas
(\ref{inj1}),(\ref{quot1}) and (\ref{inj2}),(\ref{quot2})
as an equations for the tensor product of irreducible $\slt$-modules.
To be specific,
\begin{gather}
\C^a\otimes \C^a=\C\oplus \C^{a-1}\otimes\C^{a+1} \text{ for }
(\ref{inj1}), (\ref{quot1});\\
\C^a\otimes \C^b=\C^{b-a+1}\oplus \C^{a-1}\otimes\C^{b+1}, a<b \text{ for }
(\ref{inj2}), (\ref{quot2}).
\end{gather}

The important special case of (\ref{inj1}) is the following injection:
\begin{equation}
\label{impinj}
\C^{a_1}*\ld *\C^{a_n}\hk \C^{a_1}*\ld *\C^{a_n}*\C^{a_n}*\C^{a_n}.
\end{equation}
Using (\ref{impinj}), we define an inductive limit
\begin{equation}
\label{deflim}
\C^{a_1}*\ld *\C^{a_n}*(\C^{a_n})^{*2\infty}=
\lim_{s\to\infty} \C^{a_1}*\ld *\C^{a_n}*(\C^{a_n})^{*2s}.
\end{equation}
By definition, $\C^{a_1}*\ld *\C^{a_n}*(\C^{a_n})^{*2\infty}$
is $\sltc$-module. We will show (see proposition \ref{closed}) that in fact,
it has a structure of the representation of
$\slth=\slt\otimes\C[t,t^{-1}]\oplus c\C$,
where $c$ -- a central element -- acts by the multiplication on $a_n-1$.
In order to see the connection between $\slth$-modules and the fusion
product, let us describe the latter in terms of the irreducible
representations of $\slth$.

Recall, the Lie algebra $\slth=\slt\otimes\C[t,t^{-1}]\oplus c\C$,
with the following commutator (here $(\cdot,\cdot)$ is a Killing form on
$\slt$, $x_i=x\otimes t^i$):
$$[c,x_i]=0,\ [x_i,y_j]=[x,y]_{i+j}+
c\cdot (x,y)i\delta_{i+j,0},\ x,y\in\slt.$$
Let $L_{i,k},\ i=0,\ld,k$ be an irreducible highest weight $\slth$-module
with a highest vector $v_{i,k}$ such that
$L_{i,k}=U(\slth)v_{i,k}$ and
\begin{equation}
\label{bir}
e_{\le 0}v_{i,k}=f_{<0}v_{i,k}=h_{<0}v_{i,k}=0,\ cv_{i,k}=kv_{i,k},
h_0v_{i,k}=iv_{i,k}.
\end{equation}
In addition, there exists an operator $d$ acting on $L_{i,k}$,
defined by
$[d,x_i]=ix_i,\ dv_{i,k}=0$.
It is clear that $L_{i,k}$ is bigraded by the action of the
operators $h_0,d$.
Let us denote the $h_0$-grading by $\deg_z$ and $d$-grading by $\deg_q$.

Recall that there exists a set of the extremal vectors $v(2i),\ i\in\Z$
in $L_{0,1}$, such that
$$v(0)=v_{0,1},\ e_{2i+1}v(2i)=v(2i+2),\ f_{-2i+1}v(2i)=v(2i-2).$$
$L_{1,1}$ also contains a set of the extremal vectors $v(2i+1),\ i\in\Z$:
$$v(1)=v_{1,1},\ e_{2i}v(2i-1)=v(2i+1),\ f_{-2i}v(2i+1)=v(2i-1).$$
For the positive integers $b_1\ge\ld\ge b_k$ define a subspace
$M(b_1,\ld,b_k)$ of
$L_{p_1,1}\otimes\ld\otimes L_{p_k,1}$, where $p_i$ are $0$ or $1$, and
parities of $p_i$ and $b_i$ coincide. Namely,
$$M(b_1,\ld,b_k)=
  U(\slt\otimes\C[t^{-1}])(v(-b_1)\otimes\ld\otimes v(-b_k)).$$
It turns out that $M(b_1,\ld,b_k)$ is isomorphic to some fusion product.
To make a precise statement, let us note that by definition, a fusion
product is $\slt\otimes\C[t]$-module.
But in what follows we will consider our fusions as
$\slt\otimes\C[t^{-1}]$-modules (replacing $t\to t^{-1}$). Taking this
remark into account, we can make the following statement:
\begin{utv}
\label{extr}
$$M(b_1,\ld,b_k)
 \cong (\C^2)^{*(b_1-b_2)}*(\C^3)^{*(b_2-b_3)}*\ld
 *(\C^k)^{*(b_{k-1}-b_k)}*(\C^{k+1})^{*b_k}$$
as $\slt\otimes\C[t^{-1}]$-modules.
\end{utv}

For $A=(a_1,\ld,a_n)$ introduce an associated set
$$D=(d_1,\ld,d_{k+1})\in\N^{k+1}:\
d_i=\#\{j:\ a_j=i\}, k+1=a_n=\max\{a_i\}.$$ Denote the limit from
(\ref{deflim}) as $L^D$:
\begin{equation}
L^D=\lim_{s\to\infty} (\C^1)^{*d_1}*(\C^2)^{*d_2}*\ld
*(\C^{a_n-1})^{*d_{a_n-1}}*(\C^{a_n})^{*(d_{a_n}+2s)}
\end{equation}
(we can consider $d_{a_n}$ to be equal to $0$ or $1$).
We will show that $L^D$ is a submodule in the tensor product of the
irreducible level $1$ $\slth$-modules, generated by the $\slth$ action
from some special vector.
In order to prove it, we will use the realization of the irreducible
level $1$ representations in the fermionic space (in the space of
the semi-infinite forms). Let us give a brief description of the latter.

Let $\psi(i),\phi(j),i,j\!\in\!\Z$ be an anticommuting variables
(fermions). In particular, $\psi(i)^2=\phi(j)^2=0$.
The basis of $F$ form so-called "semi-infinite" monomials in the variables
$\phi(i),\psi(j)$. They are the expressions of the form
\begin{gather}
\ld\psi(N-1)\phi(N-1)\psi(N) \phi(N) \psi(i_1)\ld\psi(i_k)
\phi(j_1)\ld\phi(j_l),\\
N<i_1<\ld<i_k, N<j_1<\ld<j_l.\notag
\end{gather}
Operators $\psi(i),\phi(j)$ acts on $F$ by the multiplication on the
corresponding fermion. Also there are the differentional operators
$\psi(i)^*,\phi(j)^*$, acting on $F$ by the differentiation
with respect to the variables
$\psi(-j),\phi(-i)$ respectively. Thus, the following commutation relations
are true ($[a,b]_+=ab+ba$):
$$[\psi^*(i),\psi(j)]_+=\delta_{i+j,0},
[\phi^*(i),\phi(j)]_+=\delta_{i+j,0}.$$
The $\slth$-action on $F$ is the following (we use the notation
$a(z)=\sum_{i=-\infty}^{\infty} a_i z^i$):
$$e(z)=\psi(z)\phi(z), f(z)=\psi^*(z)\phi^*(z),
h(z)=:\phi(z)\phi^*(z):+:\psi(z)\psi^*(z):, c=\Id.$$
Define a set of vectors $v(i)\in F$:
$$v(2N)=\ld\psi(N-1)\phi(N-1)\psi(N);\
v(2N+1)=\ld\psi(N-1)\phi(N-1)\psi(N)\phi(N).$$
We will show that there is an embedding
$\C^{a_1}*\ld*\C^{a_n}\hk F^{\otimes (a_n-1)}$ .
To be precise,
\begin{equation}
\C^{a_1}*\ld*\C^{a_n}\cong U(\slt\otimes\C[t^{-1}])
(v(-d_2-\ld-d_{a_n})\otimes\ld\otimes v(-d_{a_n}))
\label{pieces}
\end{equation}
(see statement \ref{extr}).

Using this fermionic realization, one can see that for the
irreducible $\slth$-modules we have (see also \cite{fusion}):
$$L_{i,k}=\C^{i+1}*(\C^{k+1})^{*2\infty}$$
as a subspaces of $F^{\otimes k}$.
Really, it is easy to see that $L_{i,k}$ can be embedded in
$F^{\otimes k}$ as a union of the "pieces" of the form (\ref{pieces}).
Actually, for any $D$ the module $L^D$ is a subspace of $F^{\otimes k}$.

Thus, having a vector $D\in\N^{k+1}$, we construct the level $k$ $\slth$-module
$L^D$. In our paper we find a basis of
$L^D$.  Let us explain the construction  of the basis
in the case of the vacuum irreducible representation $L_{0,k}$.

Let $v_{0,k}$ be the highest vector in $L_{0,k}$.
We realize $L_{0,k}$ as a union of the cyclic
$\slt\otimes\C[t^{-1}]$-modules with a cyclic vectors
$w_s, s=1,2,\ld.$ In our special case
$$w_s=f_1^k f_3^k \ld f_{2s-1}^k v_{0,k}.$$
Denote by $e^i(j)$ the coefficient in front of $z^j$ in the series $e(z)^i$.
We prove that the union ($s=1,2,\ld$) of the following vectors form
the $L_{0,k}$ basis:
\begin{gather}
e^k(I^k)\ld e(I^1)w_s=e^k(i^k_1)\ld e^k(i^k_{l_k}) e^{k-1}(i^{k-1}_1)\ld
e^{k-1}(i^{k-1}_{l_{k-1}}) \ld e(i^1_1)\ld e(i^1_{l_1})w_s,\\
i^{\al}_j+2\al\le i^{\al}_{j+1},\notag\\
 i^{\al}_j\ge (-2s+1)\al+\al(l_k+\ld+l_{\al+1})(l_k+\ld+l_{\al+1}+1).\notag
\end{gather}

As any integrable highest weight $\slth$-module, $L^D$ can be
decomposed into the sum of the irreducible components.
Thus, one has the decomposition:
$$L^D=M_1\otimes L_{0,k}\oplus\ld\oplus M_{k+1}\otimes L_{k,k},\
\dim M_i=c_{i,D}.$$ We show that $M_i$ are graded spaces and
numbers $c_{i,D}$ can be found in terms of the Verlinde algebra
$\V_{k+1}$, associated with $\slt$. In addition, it follows from our
results that the characters of $M_i$ are the restricted Kostka polynomials
(see \cite{kkost}, \cite{skost}). Let us recall the definition of the
Verlinde algebra.

Consider an algebra with a basis $\pi_1,\pi_2,\ld$
and multiplication $(i\le j)\ \pi_i \pi_j=\pi_{j-i+1}+\ld+\pi_{i+j-1}$
(our generators multiplies as a finite-dimensional irreducible
$\slt$-modules). By definition,
$\V_k=\langle\pi_1,\pi_2,\ld\rangle/(\pi_{k+1})$. We show that in
$\V_{k+1}$ the following equation is true:
$$\pi_D=\pi_1^{d_1}\ld \pi_k^{d_k}\pi_{k+1}^{d_{k+1}}=\sum_{i=1}^{k+1}
c_{i,D}\pi_i.$$
The proof consists of the checking that the defining relations of the
Verlinde algebra hold for our modules $L^D$. To be specific, we
check that
\begin{equation}
\text{if } \pi_D=\pi_{D'}+\pi_{D''}, \text{ then } L^D=L^{D'}\oplus L^{D''}
\label{eq}
\end{equation}
for some choice of $D,D',D''$.
The set of the triples $D,D',D''$ is suggested by the formulas (\ref{quot1}),
(\ref{quot2}).

Using the character formula from \cite{mi}, we obtain the formulas for the
character of $\C^{a_1}*\ld*\C^{a_n}$ and $L^D$. In particular, this
gives us the formula for the characters of the irreducible
representations $L_{i,k}$ from \cite{sto}. Using the characters of $L^D$,
we obtain some part of the equations $L^D=L^{D'}\oplus L^{D''}$
by the only combinatorial methods.

In the end, let us recall the method of studing of fusion products from
\cite{mi}. We use the below construction as one of the main methods.

Recall that
$$\C^{a_1}*\ld*\C^{a_n}\cong U(\slt\otimes\C[t^{-1}])w.$$
One can see that all module is generated by the only action of
the operators $e_0,\ld,e_{-n+1}.$ Thus, we have
$$\C^{a_1}*\ld*\C^{a_n}\cong \C[e_{-n+1},\ld,e_0]/I^A=M^A,$$
where $I^A$ is some ideal in the ring $\C[e_{-n+1},\ld,e_0]$.
Consider the rings isomorphism:
$$opp: \C[e_{-n+1},\ld,e_0]\to\Cn,\ opp(e_i)=e_{i+n-1}.$$
Denote
$$J^A=opp(I^A), W^A=\Cn/J^A.$$ Then the following defining relations are
true in $W^A$ (here $e^{(n)}(z)=\sum_{i=0}^{n-1} e_i z^i$):
$$e^{(n)}(z)^i\div z^{\sum_{j=1}^n (i+1-a_j)_+},\ i=1,2,\ld,$$
where $s_+=0$ if $s\le 0$;\ $s_+=s$ if $s>0$; and for polynomials
$p,q$ we write $p\div q$ if $p$ is divisible on $q$.
Our condition means that the first $-1+\sum_{j=1}^n (i+1-a_j)_+$
coefficients of
$e^{(n)}(z)^i, i=1,2,\ld$, generate the ideal $J^A$.

Our work is organized in the following way:

In the second section we construct the fermionic space $F$ as $\slth$-module
(subsection $2.1$) and describe the inclusion of $W^A$ in $F^{\otimes k}$
(theorem $2.1$).

In the third section we construct $\slth$-module $L^D$ (proposition \ref{defL})
and describe its basis (theorem \ref{basisL}).

In the fourth section we prove the decomposition formula
$L^D=L^{D'}\oplus L^{D''}$ (theorem \ref{decfun}),
and also establish the connection between our modules and Verlinde
algebra (proposition \ref{Ver}).

In the last section we write the character formula for $L^D$ (theorem
\ref{chL}) and
as a corollary  we obtain a formula for the characters of the irreducible
representations
$L_{i,k}$ (corollary \ref{irrch}). In addition, we prove by the combinatorial
means the equation
$L^D=L^{D'}\oplus L^{D''}$ (proposition \ref{decch}).
\begin{rem}
In the terminology of the works \cite{sup},\cite{mult} the theorem
\ref{chL} contains a
"fermionic" formula for the character of $L^D$. There also exists a
"bosonic" (alternating) formula, which is connected with a geometry of the
flag manifolds.
We hope to return to this topic later.
\end{rem}
\noindent Acknowledgements. The first author was partially supported by the
following grants:
CRDF RP1-2254,\ INTAS 00-55,\ RFBR 00-15-96579.

\section{Fermionic picture}
\subsection{Fermionic space (the space of the semi-infinite forms).\\}

Consider two set of variables $\psi(i), \phi(j),\ i,j\!\in\!\Z$, such
that for all $i,j\in\Z$ we have
\begin{gather}
\label{anticom}
  \psi(i)\psi(j)=-\psi(j)\psi(i),\
  \phi(i)\phi(j)=-\phi(j)\phi(i),\\
  \psi(i)\phi(j)=-\phi(j)\psi(i),\
  \psi(i)\psi(i)=\phi(i)\phi(i)=0.\notag
\end{gather}

Let
$F$ be a linear space with a basis
\begin{gather}
\ld\psi(N-2)\phi(N-2)\psi(N-1)\phi(N-1)\psi(N)\phi(N)
              \psi(i_1)\ld\psi(i_s)\phi(j_1)\ld\phi(j_t),\\
s,t\in\N,\ N,i_\al,j_\be\in\Z,\ N<i_1<\ld<i_s,\ N<j_1<,\ld<j_t.\notag
\end{gather}

In other words, we consider the space, spanned by such infinite monomials
in the variables $\phi(i),\psi(i)$
that there exists an integer $N$,
such that all the $\phi(\le N), \psi(\le N)$
are factors in it.
In addition, numbers of all factors are less than some natural number.
Let us call  the monomials with the above properties configurations,
and let the product $\prod_{i=-\infty}^{N} \psi(i)\phi(i)$
be a tale of the configuration.
It is clear that the fermions, which are not in the tale, can
be reordered, using the commutaion relations.

Define an action of the operators $\psi(i), \phi(i), \psi^*(i),\phi^*(i)$
on $F$.
By the definition, $\psi(i)(v)\ (\phi(i)(v))$, where $v$ is a
configuration, is a monomial $v\psi(i)\ (v\phi(i))$
(to get a basis element probably one must reorder the fermions).
For example, if $v$ contains $\psi(i)\ (\phi(i))$,
then $\psi(i)v=0\ (\phi(i)v=0)$.

In its turn, the operators $\psi^*(i),\phi^*(i)$ act by
the differentiation with respect to the variables
$\psi(-i),\phi(-i)$ respectively.
The following commutation relations hold ($[a,b]_+=ab+ba$):
\begin{gather}
\label{diff}
[\psi^*(i),\psi(j)]_+=\delta_{i+j,0},\ [\phi^*(i),\phi(j)]_+=\delta_{i+j,0},\\
[\psi^*(i),\phi(j)]=[\phi^*(i),\psi(j)]=0.\notag
\end{gather}
In addition, all $\psi^*(i),\phi^*(j)$ mutually anticommute.

Introduce the generating functions:
\begin{gather}
\psi(z)=\sum_{i=-\infty}^{\infty} \psi(i)z^i,\
\phi(z)=\sum_{i=-\infty}^{\infty} \phi(i)z^i,\\
\psi^*(z)=\sum_{i=-\infty}^{\infty} \psi^*(i)z^i,\
\phi^*(z)=\sum_{i=-\infty}^{\infty} \phi^*(i)z^i.\notag
\end{gather}

Recall the construction of the
$\slth=\slt\otimes\C[z,z^{-1}]+c\C$ action on $F$. The central element
acts by $1$ ($F$ is a level $1$ representation).
Define an action of the generating functions of the elements
$e,f,h$ on $F$ as follows
(for $x\in\slt$ denote $x_i=x\otimes z^i$):
\begin{gather}
e(z)=\sum_{i=-\infty}^{\infty} e_iz^i=\psi(z)\phi(z);\
f(z)=\sum_{i=-\infty}^{\infty} f_iz^i=\psi^*(z)\phi^*(z),
\end{gather}
and the action of $h(z)$ can be defined using the commutation relations
in $\slth$.

Introduce a set of vectors $v(i)\in F$:
\begin{gather}
v(2N)=\ld\psi(N-2)\phi(N-2)\psi(N-1)\phi(N-1)\psi(N);\\
v(2N+1)=\ld\psi(N-1)\phi(N-1)\psi(N)\phi(N).\notag
\end{gather}
Note that $e_{N+1}v(N)=v(N+2), e_{\le N}v(N)=0$.

Note also that the space
$\underbrace{F\otimes\ld\otimes F}_k$ inherit a structure of level $k$
$\slth$-module.

\subsection{A fermionic realization of $\Cn$-modules.\\}
Recall that in \cite{mi} we define a set of $\Cn$-modules $W^A$,
depending on the vector $A\in \N^n$:
$W^A=\Cn/J^A$, where $J^A$ is an ideal in the ring $\Cn$,
generated by the conditions on the generating function
$e^{(n)}(z)=\sum_{i=0}^{n-1} e_iz^i$:
\begin{gather}
e^{(n)}(z)^i\div z^{\sum_{j=1}^n (i+1-a_j)_+}
\end{gather}
(we use the notation $s_+=0$, if $s<0$, and $s_+=s$, if $s\ge 0$;\
$p\div q$ means that $p$ is divisible on $q$).
In other words, the ideal $J^A$ is generated by the first
$-1+\sum_{j=1}^n (i+1-a_j)_+$ coefficients of the polynomials
$e^{(n)}(z)^i, i=1,2,\ld$.

Let $A=(\an)$ and $a_1\le\ld\le a_n$. In \cite{mi} the following statement
was proved:
\begin{utv}
$W^{a_1,\ld,a_{n-1}}\hk W^A$ and
$W^A/W^{a_1,\ld,a_{n-1}}\cong W^{a_1,\ld,a_{n-1},a_n-1}$.
\end{utv}
By induction, one can obtain that $\dim W^A=\prod_{i=1}^n a_i$.

Let $k+1=a_n=\max\{a_i:i=1,\ld,n\},\ d_j=\#\{i: a_i=j\}, j=1,\ld,k+1.$
We say that the set $D=(d_1,\ld,d_{k+1})$ is associated with $A$.
Define a vector $v_A\in F^{\otimes k}$:
\begin{gather}
v_A=v(-1+d_1)\otimes v(-1+d_1+d_2)\otimes\ld\otimes v(-1+d_1+\ld+d_k).
\end{gather}
Denote $\W^A=\Cn(v_A),\ \W^A\hk F^{\otimes k}$.
Let us prove the following inclusion theorem:
\begin{theorem}
\label{real}
$\W^A\cong W^A$ as $\Cn$-modules.
\end{theorem}
\begin{proof}
First, we prove that in $F^{\otimes k}$
$$e^{(n)}(z)^i(v_A)\div z^{\Sum_{j=1}^n (i+1-a_j)_+}$$
(in other words, first $\sum_{j=1}^n (i+1-a_j)_+-1$ coefficients of
$e^{(n)}(z)^i(v_A)$ are equal to zero).
Note that
$$e(z)^i v_A\div z^{\Sum_{j=1}^n (i+1-a_j)_+}, \text{ if } i\le k;\ \
e(z)^i v_A=0\text { if } i>k
$$
(first is a consequence of the condition $e_{\le N}v(N)=0$,
and second is true, because $e(z)^{k+1}=0$ in $F^{\otimes k}$).
So, for all $i$ we obtain
$$e(z)^i v_A\div z^{\Sum_{j=1}^n (i+1-a_j)_+}.$$
We have $(x=-\sum_{i=n}^{\infty} z^{i-n}e_i)$:
$$e^{(n)}(z)^i v_A=(e(z)+z^n x)^i v_A=
         \sum_{j=0}^i \binom{i}{j} (e(z)^{i-j} z^{nj}x^j v_A).$$
Note that
$$(e(z)^{i-j} z^{nj})v_A \div z^{\Sum_{l=1}^n (i+1-a_l)_+},$$
because
$$nj+\sum_{l=1}^n (i-j+1-a_l)_+ \ge \sum_{l=1}^n (i+1-a_l)_+.$$
Thus, we obtain $J^A v_A=0$  and as a corollary,
$\dim \Cn v_A\le \prod_{i=1}^n a_i$.
To prove our theorem, it is enough to show that
$\dim\W^A=\prod a_i$.\\
In the following lemma we use the notations:
$$Aa_{n+1}=(a_1,\ld,a_n,a_{n+1}),\ A(a_{n+1}-1)=(a_1,\ld,a_n,a_{n+1}-1).$$
\begin{lem}
Let $A=(\an),\ a_1\le\ld\le a_n$. Let $a_{n+1}\ge a_n$. Denote
$B_1\subset \Cn,\ B_2\subset \C[e_0,\ld,e_n]$ such sets of monomials that
$B_1v_A=\{xv_A,\ x\in B_1\},
B_2 v_{A(a_{n+1}-1)}=\{xv_{A(a_{n+1}-1)},\ x\in B_2\}$ are linearly
independent sets of vectors in $\W^A$ and $\W^{A(a_{n+1}-1)}$ respectively.
Let
$e_nB_2=\{e_nx,\ x\in B_2\}$. Then
$B_1v_{Aa_{n+1}}\cup (e_nB_2) v_{Aa_{n+1}}$
is a linearly independent set of vectors in $\W^{Aa_{n+1}}$.
\end{lem}
\begin{proof}
Consider $3$ different cases:
$$a_{n+1}>a_n+1,\ a_{n+1}=a_n+1,\ a_{n+1}=a_n.$$
\underline{1.\ $a_{n+1}>a_n+1$.}\\
Note that in this case
\begin{gather}
v_1=v_A=v(-1+d_1)\otimes\ld\otimes v(-1+d_1+\ld+d_{a_n-1}),\\
v_2=v_{A(a_{n+1}-1)}=v_A\otimes v(-1+n)^{\otimes (a_{n+1}-a_n-1)},\notag\\
v_3=v_{Aa_{n+1}}=v_{A(a_{n+1}-1)}\otimes v(-1+n).\notag
\end{gather}
Let $b_2\in B_2,\ b_2=e_n^i x,\ x\in\Cn$. Recall that
\begin{gather}
e_{<n}v(-1+n)=0,\ e_nv(-1+n)=v(1+n),\ e_n^2v(-1+n)=0.
\end{gather}
But the sum of the summands in $e_n b_2 v_3$, with a property that
the last factor has been changed (i.e. it equals not $v(-1+n)$, but
$v(1+n)$), equals to
\begin{equation}
\label{B_2}
(i+1)b_2v_2\otimes v(1+n).
\end{equation}
Thus, since the set of vectors $B_2v_2$ is linearly independent, we obtain
the linearly independence of $e_nB_2v_3$.
In its turn, the linearly independence of $B_1v_3$ is obvious, because
$B_1\subset \Cn$ and so for $b_1\in B_1$
\begin{gather}
\label{B_1}
b_1v_3=b_1v_1\otimes v(-1+n)^{\otimes (a_{n+1}-a_n)}.
\end{gather}
Now, we must prove that $B_1v_3$ and $e_nB_2v_3$ are linearly independent
as a whole. But that is a consequence from the formula $(\ref{B_1})$
and speculations, leading to the formula $(\ref{B_2})$. \\
\underline{2.\ $a_{n+1}=a_n+1$}\\
In this case
\begin{gather}
v_A=v(-1+d_1)\otimes\ld \otimes v(-1+d_1+\ld+d_{a_n-1}),\\
v_{A(a_{n+1}-1)}=v_A,\notag\\
v_{Aa_{n+1}}=v_A\otimes v(-1+n).\notag
\end{gather}
The proof is the same as in the previous case.\\
\underline{3.\ $a_{n+1}=a_n$}\\
In this case
\begin{multline}
v_1=v_A=v(-1+d_1)\otimes\ld\otimes v(-1+d_1+\ld+d_{a_n-1}),\\
v_3=v_{Aa_{n+1}}=v_A, \\
v_2=v_{A(a_{n+1}-1)}=v(-1+d_1)\otimes\ld\otimes
v(-1+d_1+\ld+d_{a_n-2})\otimes\\
\otimes v(-1+d_1+\ld+d_{a_n-2}+d_{a_n-1}+1).
\end{multline}
Let us prove that $e_nB_2 v_3$ is linearly independent.
Recall that all of our vectors are the elements of the tensor powers
$F^{\otimes l}$. We denote the fermions, used in the construction of the
$s$-th factor as $\psi_s(i), \phi_s(j)$.

Consider a vector $w\in F^{\otimes (a_n-1)}$, defined as follows:
\begin{gather}
\label{case}
w=\phi_{a_n-1} (n-(d_1+\ld+d_{a_n-1})/2)v_2,
         \text{ if $d_1+\ld+d_{a_n-1}$ is even;}\\
 w=\psi_{a_n-1} (n-(d_1+\ld+d_{a_n-1}-1)/2)v_2,
\text{ if $d_1+\ld+d_{a_n-1}$ is odd.}\notag
\end{gather}
Thus, $w=\chi v_2$, where $\chi$
is defined by $(\ref{case})$.
One can see that $w$ is a summand in $e_nv_3$. Note also that
\begin{gather}
\label{w}
\C[e_0,\ld,e_n]v_2\cong \C[e_0,\ld,e_n]w.
\end{gather}
In fact, we multiply $v_2$ on the fermion with a such number that beeing
added to the maximal number of the opposite fermion ($\psi$ is opposite to
$\phi$ and $\phi$ to $\psi$), which is a multiplier in the last factor
in $v_2$, it alredy gives $n$. So, $\chi$ cannot appear, while acting
by the polynomials in the variables $e_0,\ld,e_n$ on $v_2$.

Now, let $l$ be a linear combination of the elements of $B_2$.
One can see that the summands of $e_nl\cdot v_3$, which contain $\chi$
as a multiplier in the last factor, arise only while acting $e_n$
in the last factor of the tensor power. In addition, if
$b_2\in B_2,\ b_2=e_n^i x,\ x\in\Cn$, then the sum of the summands
in $e_nb_2v_3$, which contain $\chi$ in the last factor, equals $(i+1)b_2w$.
Now, we obtain the linearly independence of $e_nB_2v_3$ as a consequence of
the linearly independence of $B_2v_2$ and equality (\ref{w}).

To finish the proof of the lemma, note that, firstly, $B_1v_3$
is linearly independent since $v_3=v_1$, and, secondly,
$B_1v_3\cup e_nB_2v_3$
is linearly independent since $B_1\subset\Cn$ and thus $\chi$ cannot appear
in $b_1v_3,\ b_1\in B_1$.
\end{proof}
We have proved the lemma. Now we obtain our theorem by induction on the sum
of $a_i$ (including the existense of the monomial basis in the inductive
assumption).
\end{proof}

\section {$\slth$-modules. Bases.}
\subsection{Construction of the $\slth$-modules.\\}
Let $A=(a_1\le\ld\le a_n), k+1=a_n$.
Consider the image
$I^A\in\C[e_{-n+1},\ld,e_0]$ of the ideal $J^A$ under the action of the
ring isomorphism
\begin{gather}
\Cn\to\C[e_{-n+1},\ld,e_0],\ e_i\mapsto e_{i-n+1}.
\end{gather}
As it was said in the introduction,
the corresponding quotient is isomorphic to the fusion product of the
$\slt$-modules. To be specific,
$$\C[e_{-n+1},\ld,e_0]/I^A\cong \C^{a_1}*\ld *\C^{a_n}
\text{ as $\slt\otimes (\C[z^{-1}]/z^{-n})$-modules}$$
(actually, $\C[e_{-n+1},\ld,e_0]/I^A$ has a structure of
$\slt\otimes (\C[z^{-1}]/z^{-n})$-module, see \cite{mi}).
From the theorem (\ref{real}) we obtain the following statement:
\begin{utv}
Let
\begin{gather}
w_A=v(d_1-n)\otimes v(d_1+d_2-n)\otimes\ld
     \otimes v(d_1+\ld+d_k-n),\\
M^A=\C[e_{-n+1},\ld,e_0]w_A.\notag
\end{gather}
Then as $\slt\otimes (\C[z^{-1}]/z^{-n})$-modules
$M^A\cong \C[e_{-n+1},\ld,e_0]/I^A$.
\end{utv}
\noindent Note that since $d_1+\ld+d_{k+1}=n$, then
$$w_A=v(-d_2-\ld\-d_{k+1})\otimes v(-d_3-\ld-d_{k+1})\otimes\ld
     \otimes v(-d_{k+1}).$$
Denote $A_s=A\underbrace{a_n,\ld,a_n}_{2s},\ A_s\in\N^{n+2s}.$
We will prove that
$M^{A_s}\hk M^{A_{s+1}}$ as a subspaces of $F^{\otimes k}$.

\begin{prop}
$w_A\in U(\slt\otimes\C[z^{-1}]) w_{A_1}$.
\label{inj}
\end{prop}

\begin{proof}
Recall that $d_j=\#\{i:\ a_i=j\}$. Introduce a notation:
$$\be_i=d_{i+1}+\ld+d_{k+1},\ i=1,\ld,k.$$
Then we have
\begin{gather}
w_{A}=v(-\be_1)\otimes\ld\otimes v(-\be_k),\\
w_{A_1}=v(-\be_1-2)\otimes\ld\otimes v(-\be_k-2).\notag
\end{gather}
Note that the vectors
$$w_A, e_1 w_A,\ld, e_1^{\sum_{i=1}^k (\be_i+1)} w_A$$
form an $\slt$-module with respect to the
$\slt=\langle e_1, h_0+k\cdot\Id, f_{-1}\rangle$.
From one hand, it can be checked by the direct calculation; from the other
one can see that there exists an isomorphism
$$\C[e_1,e_0,\ld,e_{-n+1}]w_A\cong \C[e_0,\ld,e_{-n}]/I^{Aa_n},\
e_i\mapsto e_{i-1}.$$
But in any fusion product the $e_0$-action on the highest vector
spannes the $\slt$-module.

Note that
\begin{equation}
e_1^{\sum_{i=1}^k (\be_i+1)} w_A \text{ is proportional to }
v(\be_1+2)\otimes\ld\otimes v(\be_k+2).
\end{equation}
Also, one can see that
\begin{equation}
e_0^{\sum_{i=1}^k (\be_i+2)} w_{A_1}
\text{ is proportional to } v(2+\be_1)\otimes\ld\otimes v(2+\be_k).
\end{equation}
Thus, we obtain the following equality
\begin{equation}
\const\cdot w_A=
f_{-1}^{\sum_{i=1}^k (\be_i+1)}
e_0^{\sum_{i=1}^k (\be_i+2)} w_{A_1}
\end{equation}
(of course, the constant doesn't equal to zero). Proposition is proved.
\end{proof}

\begin{rem}
\label{act}
From the proof we obtain that
$e_1 M_A\hk M_{A_1}$.
\end{rem}

\begin{rem}
\label{act2}
From the proof of the proposition one can obtain that
$$w_{A_1}\hk U(\slt\otimes\C[z])w_A.$$
\end{rem}

We have proved that $w_A\in U(\slt\otimes\C[z^{-1}]) w_{A_1}$.
Thus, we have a chain of the inclusions:
$$M^A\hk M^{A_1}\hk M^{A_2}\hk\ld.$$
\begin{opr}
\label{defL}
Let $D\in {(\N\cup 0)}^{k+1}$, and $d_{k+1}=0$ or $d_{k+1}=1$.
Let $d_1+\ld+d_{k+1}=n$. Consider $A\in\N^{n+2}$, such that
$d_j=\#\{i:\ a_i=j\}$ for $j=1,\ld,k$ and $\#\{i:\ a_i=k+1\}=d_{k+1}+2$.
Then
$$L^D=L^{d_1,\ld,d_{k+1}}=\bigcup_{s=0}^{\infty} M^{A_s}.$$
Notation:
$$L^D=
(\C^1)^{*d_1}*(\C^2)^{*d_2}*\ld*(\C^k)^{*d_k}*
(\C^{k+1})^{*(d_{k+1}+2\infty)}.$$
\end{opr}
By definitions, $L^D$ is only $\slt\otimes\C[z^{-1}]$-module.
We show that in fact it has a structure of $\slth$-module.
\begin{prop}
\label{closed}
$L^D$ is closed with repsect to the action of $\slth$.
\end{prop}
\begin{proof}
By the remark (\ref{act}), $L^D$ is closed with respect to
the action of $e_1$.
In addition, $L^D$ is $\slt\otimes\C[z^{-1}]$-module. But $\slth$
is generated by
$\slt\otimes\C[z^{-1}]$ and $e_1$. Proposition is proved.
\end{proof}
It is clear that $L^{d_1,\ld,d_{k+1}}$ is an integrable level $k$
representation. It is cyclic $\slth$-module and may be generated from any of
the vectors
$w_{A_s}$ (see the proposition \ref{inj} and remark \ref{act2}).

\begin{cor}
Let $0\le i\le k$. Define such $D$ that  $d_j=\delta_{i+1,j}$.
Then $L^D\cong L_{i,k}$.
We can write it as $L_{i,k}=\C^{i+1}*(\C^{k+1})^{*2\infty}$.
\end{cor}
\begin{rem}
One can see that $L^{d_1,d_2,\ld,d_{k+1}}=L^{d_2,\ld,d_{k+1}}$.
That means that the number of the one-dimensional representations is
not important.
\end{rem}
\begin{rem}
Let $p(D)$ be a number of the odd numbers in the sequence
$$d_{k+1}, d_{k+1}+d_k,\ld, d_{k+1}+\ld+d_2.$$
Then $$L^D\hk L_{1,1}^{\otimes p(D)}\otimes L_{0,1}^{\otimes (k-p(D))}.$$
\end{rem}

\subsection{Basis of $L^D$.\\}
Here we construct an $L^D$ basis.
First we will construct the basis of the space
$W^{A_\infty}=\C[\ld,e_{-1},e_0,e_1,\ld]v_A$.

Let $k+1=a_n=\max\{a_i\}$. Recall that
$e(z)=\sum_{i=-\infty}^{\infty} e_i z^i.$
Denote by $e^i(j)$ the coefficient in $e(z)^i$ in front of $z^j$.
Note that $e_{<0}v_A=0$,
so $e(z)v_A=(\sum_{i=0}^{\infty} e_i z^i) v_A.$
Let us prove the following lemma:
\begin{lem}
The following elements are linearly independent in $W^{A_\infty}$:
\begin{gather}
\label{basis}
e(I^1)\ld e^k(I^k)v_A=e(i^1_1)\ld e(i^1_{l_1}) e^2(i^2_1)\ld
e^2(i^2_{l_2}) \ld e^k(i^k_1)\ld e^k(i^k_{l_k})v_A,\\
\label{condition}
i^{\al}_j+2\al\le i^{\al}_{j+1},\\
 i^{\al}_j\ge \al d_1+(\al-1)d_2+\ld+d_\al+
 \al(l_k+\ld+l_{\al+1})(l_k+\ld+l_{\al+1}+1).\notag
\end{gather}
\end{lem}
\begin{proof}
Recall that $v_A\in F^{\otimes k}$. We denote the fermions, used in the
construction of the $s$-th factor of $F^{\otimes k}$ as
$\psi_s(i),\phi_s(j)$.
Denote $e_{i,s}$ the following operators in $F^{\otimes k}$:
$$e_{i,s}=\underbrace{\Id\otimes\ld\otimes\Id}_{s-1}\otimes e_i\otimes\Id
\otimes\ld\otimes\Id.$$
Let $[x]$ be an integral part of $x$.
Define an operator
$$m(i^{\al}_1,\ld,i^{\al}_{l_\al},d_1,\ld,d_{\al})=
m(I^{\al},d_1,\ld,d_{\al}):F^{\otimes k}\to F^{\otimes k}$$
by the following way:
\begin{multline}
\prod_{j=1}^{\al-1}\prod_{m=1}^{l_{\al}} e_{d_1+\ld+d_j+2(m-1),j}
\times\\
\times\prod_{m=1}^{l_{\al}}
\psi_{\al}\left(\left[\frac{i^{\al}_m-(\al-1)d_1-(\al-2)d_2-\ld-d_{\al-1}-
(\al-1)2(m-1)+1}{2}\right]\right)\times\\
\times\prod_{m=1}^{l_{\al}}
\phi_{\al}\left(\left[\frac{i^{\al}_m-(\al-1)d_1-(\al-2)d_2-\ld-d_{\al-1}-
(\al-1)2(m-1)}{2}\right]\right).
\end{multline}

Let $m(I^1,\ld,I^k)v_A\in F^{\otimes k}$ be the following vector:
\begin{multline}
m(I^1,\ld,I^k)v_A=m(I^1,d_1+2(l_k+\ld+l_2))
m(I^2,d_1+2(l_k+\ld+l_3), d_2)\ld\\
\ld m(I^{k-1},d_1+2l_k,d_2,\ld,d_{k-1}) m(I^k,d_1,\ld,d_k)v_A.$$
\end{multline}
Note that if $I^1,\ld,I^k$ satisfy the condition (\ref{condition}), then
$m(I^1,\ld,I^k)v_A$ is one of the summands in $e(I^1)\ld e^k(I^k)v_A$.
Let us explain why.

Let us act in a following way: first by $e^k(I^k)$, then by
$e^{k-1}(I^{k-1})$, and so on till to $e(I^1)$.
All this operators act on $F^{\otimes k}$. Every time, while acting
by $e^{\al}(I^{\al})$ on $F^{\otimes k}$ we consider only the summands
in which the action is happening on the first $\al$ factors.
In addition, the action on the first $\al~-1$ factors is happening
by the "minimal" possible way, that is we are multiplying on the fermions
with a minimal possible numbers.
In the same time the action on the $\al$-th factor is happening by the
"middle" way, that is $e^{\al}(j)$ acts on the $\al$-th factor
by the multiplication on
$\psi_{\al} [\frac{j+1}{2}]\phi_{\al} [\frac{j}{2}]$. Condition
(\ref{condition}) garanties that while acting by the described way, we will
get a nonzero vector.

Now, let us prove the linearly independence of our vectors.
Consider the linear combination
\begin{equation}
\sum \be_{I^1,\ld,I^k} e(I^1)\ld e^k(I^k)v_A,
\label{be}
\end{equation}
where all the sets $(I^1,\ld,I^k)$ satisfy the condition (\ref{condition}).
Between the sets from (\ref{be}) choose the one $(I^1_0,\ld,I^k_0)$ with
a property that for all other sets from (\ref{be}), we have:
there exists $j, m$, such that $(i_0)^j_m>i^j_m$ and for all such $j_1,m_1$,
that
$j_1>j$ or $j_1=j, m_1<m$ we have $(i_0)^j_m=i^{j_1}_{m_1}.$
Note that from the construction, we obtain that
$m(I_0^1,\ld,I_0^k)v_A$ appears in
(\ref{be}) only from $e(I_0^1)\ld e^k(I_0^k)v_A$ (and with a nonzero
coefficient) and doesn't appear from any other summand. Thus we obtain
that
$$\sum \be_{I^1,\ld,I^k} e(I^1)\ld e^k(I^k)v_A\ne 0.$$
Lemma is proved.
\end{proof}

\begin{prop}
Vectors (\ref{basis}) with a condition (\ref{condition}),
form a base of $W^{A_{\infty}}$.
\end{prop}
\begin{proof}
Since we have checked the linearly independence, it is enough
to compare the characters of the probable basis and the space
$W^{A_{\infty}}$.
One can see that the character of the probable basis is
\begin{equation}
\sum_{i_1,\ld,i_k=0}^{\infty}
         (zq^{-1})^{\sum\limits_{l=1}^k li_l}
\frac{q^{\sum\limits_{s,t=1}^k \min(s,t)i_s i_t+
\Sum_{l=1}^k d_l(i_l+2i_{l+1}+\ld+(k-l+1)i_k)}}
{\qf{i_1}\ld\qf{i_k}}.
\label{ch}
\end{equation}
Thus, the character of $W^{A_\infty}$ is more or equal
(in every homogeneous component)
than the expression (\ref{ch}). We will show that in fact we have an
equality.

Consider a quotient $B^A$ of the polynomial ring
$$\C[b_1(0),b_1(1),\ld,b_2(0),b_2(1),\ld,b_k(0),b_k(1),\ld].$$
The defining relations in $B^A$ are the following relations on the generating
functions $b_i(z)=\sum_{j=0}^\infty b_i(j)z^j$ (we use the notation
$b^{(s)}(z)$ for the $s$-th derivative of the series):
\begin{gather}
\label{dcond}
b_i(z)\div z^{id_1+(i-1)d_2+\ld+d_i},\\
b_i^{(s)}(z)b_j^{(t)}(z)=0 \text{ if } i<j,\ s+t<2i.\notag
\end{gather}

\begin{utv}
Define $\deg_z b_i(j)=i,\ \deg_q b_i(j)=j$. Then the character of $B^A$
is coinciding with the expression (\ref{ch}).
\end{utv}
\begin{proof}
It is known, (see for example \cite{coinv2}) that the dual space
$(B^A)^*$ can be realized as a sum of the spaces of polynomials
$$f(z_1^1,\ld,z^1_{i_1};\ld;z^k_1,\ld,z^k_{i_k}),\ i_\al\ge 0,$$
symmetric with respect to all of the groups of variables
$z^s_1,\ld,z^s_{i_s}$ and which can be written in a following form:
\begin{gather}
\label{pol}
f=g\prod_{l=1}^k \prod_{i=1}^{i_l} (z^l_i)^{ld_1+(l-1)d_2+\ld+d_l}
\prod_{1\le s\le t\le k}
\prod_{\gf{i=1,\ld,i_s}{j=1,\ld,i_t}} (z^s_i-z^t_j)^{2s},
\end{gather}
where $g(z_1^1,\ld,z^1_{i_1};\ld;z^k_1,\ld,z^k_{i_k})$
is an arbitrary polynomial, symmetric  with respect to all of the groups of
variables. In other words, $f$ is vanishing $2s$ times on the diagonals
$z^s_i=z^t_j$, $s\le t$.

One can see that the character of the polynomials (\ref{pol}) coincides
with (\ref{ch}).
\end{proof}

To prove the proposition, it is enough to show that
$\ch W^{A_\infty}\le \ch B^A$ (we mean that the inequalities are true
in each homogeneous component). In order to do it, we will use the
construction
for the filtration on the dual space $(W^{A_\infty})^*$, which is explained
in \cite{coinv2} in a more general case.

One can see that a dual space $(W^{A_\infty})^*$ can be realized
as a sum of the spaces of the polynomials
$f(z_1,\ld,z_s),\ s=0,1,2,\ld$, which satisfy the following conditions:
\begin{gather}
f(\underbrace{z,\ld,z}_i,z_{i+1},\ld,z_s)\div
z^{id_1+(i-1)d_2+\ld+d_i},\ i=1,\ld,k,\\
f(\underbrace{z,\ld,z}_{k+1},z_{k+2},\ld,z_s)=0.\notag
\end{gather}
(This realization is tantamount to the fact that
$W^{A_\infty}$ is a quotient of the polynomial ring $\C[e_0,e_1,\ld]$,
with a defining relations:
\begin{equation}
e^i(z)\div z^{id_1+\ld+d_i},\ i=1,\ld,k,\ e^{k+1}(z)=0.
\end{equation}
In \cite{mi} we have proved it in the case $d_1=\ld=d_k=0$. In the general case the
proof is the same.)

As it is shown in \cite{coinv2}, there is a filtration on the space
$(W^{A_\infty})^*$ by the subspaces, enumerated by the Young diagrams,
and the adjoint quotients are isomorphic to the spaces of polynomials
(\ref{pol}).
This gives us that the character of $W^{A_\infty}$ is equal to
(\ref{ch}).
\end{proof}

\begin{cor}
Elements
\begin{multline}
e(I^1)\ld e^k(I^k)w_A=e(i^1_1)\ld e(i^1_{l_1}) e^2(i^2_1)\ld
e^k(i^2_{l_2}) \ld e^k(i^k_1)\ld e^k(i^k_{l_k})w_A,\\
i^{\al}_j+2\al\le i^{\al}_{j+1},\\
 i^{\al}_j\ge
 \al-\al (d_{k+1}+\ld+d_{\al+1})-(\al-1)d_{\al}-\ld-d_2+
 (l_k+\ld+l_{\al+1})(l_k+\ld+l_{\al+1}+1)
\end{multline}
form a base of $\C[\ld,e_{-1},e_0,e_1,\ld]w_A$.
(Denote this set of vectors as $B_A$).
\end{cor}
\begin{lem}
\label{injbase}
$B_{A_s}\subset B_{A_{s+1}}$.
\end{lem}
\begin{proof}
Note that in $F^{\otimes k}$ the following is true
$$e(z)^k=k!\cdot e(z)\otimes\ld\otimes e(z).$$
Thus we obtain the following equality ($\al_i=d_2+\ld+d_i$):
\begin{equation}
e^k\left(\sum_{i=2}^{k+1} (-\al_i-1)\right) w_{A_2}=k!\cdot w_A.
\end{equation}
As a consequence, we obtain that $B_{A_s}\subset B_{A_{s+1}}$.
\end{proof}

\begin{theorem}
\label{basisL}
$\bigcup_{s=0}^{\infty} B_{A_s}$ is a basis of $L^D$.
\end{theorem}
\begin{proof}
By definition
$$L^D=\bigcup_{s=0}^{\infty} \C[\ld,e_{-1},e_0,e_1,\ld]w_{A_s}.$$
Using the lemma \ref{injbase}, we obtain the theorem.
\end{proof}

\begin{rem}
Using the theorem \ref{basisL} and formula (\ref{ch}) one can obtain the
character of $L^D$ (note that in \cite{sto} character of
$L_{i,k}$ was obtained by this method). But in the last section we will
obtain the character of $L^D$ using a fact that $L^D$ is a limit of
finite-dimensional fusion products.
\end{rem}

\section{The decomposition of $L^D$. Verlinde algebra.}
\subsection{Functional recurrent formula.\\}
Let $A=(a_1\le\ld\le a_n)\in {\N^n},\ k+1=a_n$,\
$D\in (\N\cup 0)^{k+1}$ is associated with $A$.
We will construct such $A',A''$ that
$W^{A'}\hk W^A$ and $W^A/W^{A'}\cong W^{A''}$.

Recall that $e^{(n)}(z)=\sum_{i=0}^{n-1} e_iz^i$. Denote by
$[e^{(n)}(z)^i]_j$ the coefficient in front of $z^j$ in $e^{(n)}(z)^i$.
\begin{lem}
Let $d_l\ne 0$. Then
$$[e^{(n)}(z)^{l-1}]_{\sum_{j=1}^{l-1} (l-j)d_j}v_A=
(l-1)!(e_{d_1}\otimes\ld\otimes e_{d_1+\ld+d_{l-1}}
\otimes\Id\otimes\ld\otimes\Id) v_A.$$
\label{subp}
\end{lem}
\begin{proof}
Firstly, note that
\begin{equation}
[e^{(n)}(z)^{l-1}]_{\sum_{j=1}^{l-1} (l-j)d_j}v_A=
e^{l-1}\left(\sum_{j=1}^{l-1} (l-j)d_j\right) v_A
\label{reduct}
\end{equation}
(recall that $e^i(j)$ is a coefficient in front of $z^j$ in the series
$e(z)^i$).
In fact, if $c=\sum_{i=n}^{\infty} e_iz^{i-n}$, then
\begin{equation}
e(z)^{l-1}v_A=e^{(n)}(z)^{l-1}v_A+\sum_{i=0}^{l-2}
\binom{l-1}{i} z^{n(l-1-i)}c^{l-1-i}e^{(n-1)}(z)^i v_A.
\label{pow}
\end{equation}
Recall that
$$e^{(n)}(z)^i v_A\div z^{\sum_{j=1}^i d_j(i+1-j)}.$$
We have:
$$z^{n(l-1-i)}e^{(n-1)}(z)^i v_A\div z^{n(l-1-i)+\sum_{j=1}^i d_j(i+1-j)}.$$
In the same time, since $\sum_{j=1}^{k+1} d_j=n$ and $d_l\ne 0$, then
$$n(l-1-i)+\sum_{j=1}^i d_j(i+1-j)>\sum_{j=1}^{l-1} (l-j)d_j.$$
From this and (\ref{pow}), we obtain (\ref{reduct}), because all the summands
in the right hand side of (\ref{pow}), but $e^{(n)}(z)^{l-1}$ are
divisable on $z^{1+\sum_{j=1}^{l-1} d_j(l-j)}$.

Consider now $e(z)^{l-1} v_A$. Since $e(z)^2=0$ in $F$, we have
\begin{gather}
e(z)^{l-1}v_A=(l-1)!\sum_{1\le i_1<\ld<i_{l-1}\le k}
e(z)_{i_1}\ld e(z)_{i_{l-1}}v_A,\\
e(z)_i=\underbrace{\Id\otimes\ld\otimes\Id}_{i-1}\otimes e(z)\otimes\Id
\otimes\ld\otimes\Id\notag.
\end{gather}
Consequently, since $d_l\ne 0$ and $e_{\le d} v(d)=0$, then
$$
[e(z)^{l-1}]_{\sum_{j=1}^{l-1} (l-j)d_j}v_A=
(l-1)!(e_{d_1}\otimes\ld\otimes e_{d_1+\ld+d_{l-1}}
\otimes\Id\otimes\ld\otimes\Id) v_A.$$
Lemma is proved.
\end{proof}
\begin{rem}
Note that if $l-1=k$ then we as a corollary obtain the proposition
(\ref{inj}).
\end{rem}
\begin{rem}
One can see that
\begin{multline}
(e_{d_1}\otimes\ld\otimes e_{d_1+\ld+d_{l-1}}
\otimes\Id\otimes\ld\otimes\Id) v_A
=v(1+d_1)\otimes v(1+d_1+d_2)\otimes\ld\\
\ld\otimes v(1+d_1+d_2+\ld+d_{l-1})\otimes v(-1+d_1+\ld+d_l)\otimes\ld
\otimes v(-1+d_1+\ld+d_k).
\end{multline}
Denote this vector by $u$.
\end{rem}

Using the fermionic realization of the modules $W^A$, we obtain:
\begin{utv}
\label{injfus}
a).\ Let $d_l\ge 2$. Denote by $A'\in \N^{n-2}$ such set that
the associated set $D'$ is the following: $d'_l=d_l-2,\
d'_i=d_i$ if $i\ne l$. Then, if we will consider $W^{A'}$ as a ring
$\C[e_0,\ld,e_{n-3}]/J^{A'}$, then there exists the isomorphism of the
rings:
$$\C[e_2,\ld,e_{n-1}]u\cong W^{A'},\ e_i\mapsto e_{i-2}.$$
b).\ Let $d_l=1,\ l=\min\{m:\ d_m\ne 0\},\ l\ne k+1$.
Let $d_{l+1}=\ld=d_{l_1-1}=0, d_{l_1}\ne 0$.
Denote by $A'\in \N^{n-1}$ such set that an associated set $D'$ is of
the following form:
$$d'_l=d_l-1, d'_{l_1}=d_{l_1}-1, d'_{l_1-l+1}=d_{l_1-l+1}+1,$$
and for all other $i$\ $d_i=d'_i$.
Then we have the isomorphism of the rings:
$$\C[e_1,\ld,e_{n-1}]u\cong W^{A'},\ e_i\mapsto e_{i-1}.$$
\end{utv}
\begin{rem}
Note that the condition $l=\min\{m:\ d_m\ne 0\}$ from $b)$ seems to be
unnatural.
But in the case
of geneneral $l$ the module $\C[e_1,\ld,e_{n-1}]u$ is not of the type
$W^{A'}$ for any $A'$. Still, the structure of modules
$\C[e_1,\ld,e_{n-1}]u$ in the general case is very interesting too and
we hope to study it in the future.
\end{rem}

So, there is a subspace $W^{A'}$ in $W^A$. Let us study the quotient.
Recall that (see \cite{mi}) the dual space $(W^A)^*=(\Cn/J^A)^*$
can be realized as a space of the polynomials
$f(z_1,\ld,z_s),\ s=0,1,2,\ld$ with a following conditions:\\
$1).\ f(z_1,\ld,z_s) \text{ is symmetric.}\\
2).\ deg_{z_i} f(z_1,\ld,z_s)<n.\\
3).\ f(\underbrace{z,\ld,z}_i,z_{i+1},\ld,z_s)\div
z^{\Sum_{j=1}^n (i+1-a_j)_+},\ i=1,\ld,s.$\\
Let us prove the following proposition:
\begin{prop}
\label{razl}
a).\ $d_l\ge 2$. Denote by $A''\in\N^n$ such set that
the associated set $D''$ is the following:
$$d''_{l-1}=d_{l-1}+1, d''_{l}=d_l-2, d''_{l+1}=d_{l+1}+1,
d''_{i}= d_{i} \text{ otherwise}.$$
Then $W^A/W^{A'}\cong W^{A''}$.\\
b).\ $d_l=1,\ l=\min\{m:\ d_m\ne 0\}$. Denote by $A''\in\N^n$ such set that
the associated set
$D''$ is the following:
$$
d''_{l-1}=d_{l-1}+1, d''_{l}=0,
d''_{l_1}=d_{l_1}-1,
d''_{l_1+1}=d''_{l_1+1}+1, d''_{i}= d_i \text{ otherwise}.
$$
Then $W^A/W^{A'}\cong W^{A''}$.
\end{prop}
\begin{proof}
Note that the dual space $$\left(W^A/\Cn u\right)^*$$
can be realized as a space of the symmetric polynomials $f$ in $s$ variables
($s=0,1,\ld$), with a degree not more than $n-1$ in each and satisfying
the following conditions:\\
$*).\ f(\underbrace{z,\ld,z}_i,z_{i+1},\ld,z_s)\div
z^{\Sum_{j=1}^n (i+1-a_j)_+},\ i=1,\ld,s.\\
**).\ f(\underbrace{z,\ld,z}_{l-1},z_l,\ld,z_s)\div
z^{1+\Sum_{j=1}^n (l-a_j)_+}\\$
(recall that $u=[e^{(n)}(z)]_{\sum_{j=1}^{l-1} (l-j)d_j} v_A$).
We will prove that if the polynomial $f$ satisfies the conditions $*),**)$,
then it also satisfies the conditions on the function from the dual space
$(W^{A''})^*$. If we also check that
\begin{equation}
\dim W^A=\dim W^{A'}+\dim W^{A''},
\label{dim}
\end{equation}
we will obtain the proof of the proposition.

Firstly, let us check that $\dim W^A=\dim W^{A'}+\dim W^{A''}$.
As we have mentioned above, $\dim W^A=\Prod a_i$. But then $(\ref{dim})$
is tantamount to the equality $l^2=1+(l-1)(l+1)$ in the case $a)$ and
$ll_1=l_1-l+1+(l-1)(l_1+1)$ in the case $b)$.

Now, let us prove that if $f$ satisfies the conditions $*),**)$, then
$f\in (W^{A''})^*$ that is
\begin{gather}
f(\underbrace{z,\ld,z}_i,z_{i+1},\ld,z_s)\div
z^{\Sum_{j=1}^i d''_j(i+1-j)_+},\ i=1,\ld,s.
\label{cond}
\end{gather}
$a).$\  $d_l\ge 2$.
If $i<l-1\  (\ref{cond})$ is a consequence of $*)$. If $i=l-1$ it is a
consequence from $**)$.
Let $i=l$. Condition $*)$ gives us
$$f(\underbrace{z,\ld,z}_l,z_{l+1},\ld,z_s)\div
z^{\sum_{j=1}^l (l+1-j)d_j}.$$
But $d''_{l-1}=d_{l-1}+1, d''_l=d_l-2$. So
$$2d''_{l-1}+d''_l=2d_{l-1}+d_l.$$ Thus, $(\ref{cond})$ is true for
$i=l$. Now let $i>l$. Then
\begin{gather*}
(d_{l-1}+1)(i+1-(l-1))+(d_l-2)(i+1-l)+(d_{l+1}+1)(i+1-(l+1))=\\
=d_{l-1}(i+1-(l-1))+d_l(i+1-l)+d_{l+1}(i+1-(l+1)).\notag
\end{gather*}
So, we obtain (\ref{cond}) as a corollary from $*)$.
Case $a)$ is proved.\\
$b).\ d_l=1$. Recall that for $m<l$ we have $d_m=0$. So for
$i<l-1$ the condition (\ref{cond}) is trivial.

We know from $**)$ that
\begin{equation}
f(\underbrace{z,\ld,z}_{l-1},z_l,\ld,z_s)\div z.
\end{equation}
Using this formula, one can obtain that
\begin{equation}
f(\underbrace{z,\ld,z}_{l-1+k},z_{l+k},\ld,z_s)\div z^{k+1},\ k\ge 0.
\end{equation}
This gives us (\ref{cond}) for $i=l-1,\ld,l_1-1$. In the same time, for
$i\ge l_1$ the condition (\ref{cond}) coincides with a condition $*)$.
Proposition is proved.
\end{proof}

\subsection{Decomposition of $L^D$.\\}
Here we obtain the decomposition of $L^D$ on the irreducible components.

Recall that we construct $L^D$ as a limit of the inclosed
spaces $M^{A_s}$. Note that the proposition (\ref{razl}) can be reformulated
in a following way: $M^A\hookleftarrow M^{A'},\ M^A/M^{A'}\cong M^{A''}$.

Recall that we have constructed a fermionic realization of $M^A$.
We want to understand, how do $M^{A}$ and $M^{A'}$ are embedded in
$F^{\otimes k}$ with respect to each other.
We know that $M^A/M^{A'}\cong M^{A''}$. In the same time, we will calculate
the minimal level, in which $M^A$ doesn't coincide with $M^{A'}$ that is the
"starting" level for $M^{A''}$. It is clear that the number of this level is
$\deg_q w_A-\deg_q w_{A''}$.

Let $2\al_i=d_{k+1}+\ld+d_{i+1}, i=1,\ld,k$, and define $p(D)$
as a number of half-integers $\al_i$.
\begin{lem}
\label{deg}
$a).\ d_l\ge 2, l\ne k+1.$ Then
$$\deg_q w_A-\deg_q w_{A''}=\frac{p(D'')-p(D)}{4}+\frac{d_l-1}{2}.$$
$b).\ d_l=1,\ l_1=\min\{i>l:\ d_i\ne 0\}.$ Then
$$\deg_q w_A-\deg_q w_{A''}=\frac{p(D'')-p(D)}{4}+\frac{d_{l_1}}{2}.$$
\end{lem}
\begin{proof}
From the fermionic realization one can obtain
$$deg_q w_A=\frac{-p(D)}{4}+\sum_{i=1}^k \al_i^2.$$
The rest of the proof is straightforward.
\end{proof}

\begin{theorem}
\label{decfun}
Let $D\in (N\cup 0)^{k+1}, d_1=0,d_{k+1}\le 1.$ Let
$d_2+\ld+d_{k+1}\ge 2$. Let $2\le l\le k,\ d_{<l}=0,\
l_1=\min\{i:\ i>l, d_i\ne 0\}$.\\
$a).\ d_l\ge 2$. Define $D',D''$ as above with an only change:
if $l=k$, then $d''_{k+1}=(d_{k+1}+1)\mod 2$. Then
$$\ch L^D=\ch L^{D'}+q^{\frac{p(D'')-p(D)}{4}+\frac{d_l-1}{2}}\ch L^{D''}.$$
$b_1).\ d_l=1,\ l_1\ne k+1$. Then
$$\ch L^D=\ch L^{D'}+q^{\frac{p(D'')-p(D)}{4}+\frac{d_{l_1}}{2}}\ch L^{D''}.$$
$b_2).\ d_l=1,\ l_1=k+1.$ Then
$$\ch L^D=\ch L^{D'}.$$
\end{theorem}
\begin{proof}
Our theorem is a consequence of the lemma (\ref{deg}).
Note only that in the case $b_2)$ the second summand vanishes, because
the difference
$\deg_q w_{A_s}-\deg_q w_{A''_s}$ goes to infinity while
$s\to\infty$. In the other cases this difference doesn't depend on~ $s$.
\end{proof}

\begin{cor}
From the equations for the characters of the integrable highest weight
$\slth$-modules, we obtain the equalities for the modules itself:
(all notations are from the theorem)\\
$a).\ L^D=L^{D'}\oplus L^{D''}.\
b_1).\ L^D=L^{D'}\oplus L^{D''}.\
b_2).\ L^D=L^{D'}$.
\label{decomp}
\end{cor}

\begin{cor}
\label{alg}
We obtain an algorithm for the decomposition of the modules $L^D$ on the
irredicible components.
\end{cor}
\begin{proof}
Note that such $D$ that $\sum d_i\le 1$, correspond to the irreducible
representations $L_{i,k}$.
To all the rest of $L^D$ we can apply our decomposing procedure.
Note that for $D,D',D''$, defined above, we have
$\prod i^{d_i}>\prod i^{d'_i},\ \prod i^{d_i}>\prod i^{d''_i}.$
So, our algorithm in a finitely many steps will give us
the decomposition of $L^D$ on the sum of the irreducible $\slth$-modules.
\end{proof}

\subsection{The connection with the Verlinde algebra.\\ }
Recall the definition of the Verlinde algebra $\V_k$, associated with
a Lie algebra $\slt$. Consider an algebra with a basis $\pi_1,\pi_2,\ld$
and with a multiplication
$(i\le j)\ \pi_i \pi_j=\pi_{j-i+1}+\ld+\pi_{i+j-1}$
(generators are multiplying as a finite-dimensional irreducible
$\slt$-modules). By definition,
$\V_k=\langle\pi_1,\pi_2,\ld\rangle/(\pi_{k+1})$.
Here is a list of some properties of $\V_k$:
\begin{utv} The following equations are true in $V_k$:\\
\label{verprop}
$1).\ \pi_k^2=1(=\pi_1).\\
2).\ \pi_i^2=1+\pi_{i-1}\pi_{i+1},\ i=2,\ld,k-1.\\
3).\ \pi_i \pi_j=\pi_{j-i+1}+\pi_{i-1}\pi_{j+1},\ i<j, i=2,\ld,k-2,
j=3,\ld,k-1.\\
4).\ \pi_i \pi_k=\pi_{k-i+1},\ i=2,\ld,k-1.$
\end{utv}

Let $D\in (\N\cup 0)^{k+1}, d_{k+1}\le 1$. Denote
\begin{equation}
\pi_D={\pi_1}^{d_1} {\pi_2}^{d_2}\ld {\pi_{k+1}}^{d_{k+1}}\in \V_{k+1}.
\end{equation}
We have
\begin{equation}
\pi_D=c_{1,D}\pi_1+\ld+c_{k+1,D}\pi_{k+1}.
\end{equation}
From the corollaries \ref{decomp}, \ref{alg} and statement \ref{verprop}
we obtain the following proposition:
\begin{prop}
There is a decomposition
\label{Ver}
$$L^D=M_1\otimes L_{0,k}\oplus\ld\oplus M_{k+1}\otimes L_{k,k},$$
where $\dim M_i=c_{i,D}$ and $M_i$ are graded spaces. The character of
$M_i$ can be obtained by the procedure, described above.
\end{prop}

\section{Combinatorial calculations}
In this section we will obtain the formulas for the character of $L^D$.
As a consequence we will obtain the formula from \cite{sto} for the
characters of the irreducible representations $L_{i,k}$.
Using this formulas, we will reprove an equation $a)$ from the
corollary (\ref{decomp}).

\subsection{The character formula for $L^D$.\\}
One can see that the spaces $L^D$ are bigraded by the operators
$h_0$ and $d$ (for $x\in\slt$\ $[d,x_i]=ix_i$).
So, one can introduce a character as a trace $\Tr q^d z^{h_0}$.
In the same time, it is clear that for all $D$ the eigenvalues of
$h_0$ on $L^D$ are of the same parity.
So, we will define the character in the following way:
$\ch L^D=\Tr q^d z^{\frac{h_0}{2}}$. In this definition, all the powers
of $z$ are integer or half-integer.
(One can see that they are integer, if $\sum_{i=2}^{k+1} (i-1)d_i$ is even,
and half-integer, if it is odd).

Recall that in \cite{mi} we obtain the following formula for the character
of $W^A$\ ($A=(a_1\le\ld\le a_n),\ a_n=k+1, d_j=\#\{i:\ a_i=j\}, d_1=0$):
\begin{multline}
\label{chW}
\ch(W^A,q,qz)=\\
=\sum_{{j_k}=0}^{d_{k+1}}\sum_{j_{k-1}=0}^{d_k+j_k}\ld
\sum_{j_1=0}^{d_2+j_2}
         z^{\sum\limits_{l=1}^k j_l}
q^{\sum\limits_{l=1}^k j_l(d_2+\cdots +d_l+j_l)}\times\\
\times\qb{d_{k+1}}{j_k}\qb{d_k+j_k}{j_{k-1}}\ld\qb{d_2+j_2}{j_1}.
\end{multline}
Let us now obtain the character of $M^A$ (again, powers of $z$ will be
integer or half-integer).
Denote
$$2\al_i=d_{k+1}+\ld+d_{i+1}\ i=1,\ld,k.$$
In this notations $p(D)=\#\{i:\ \al_i\in(\Z+\frac{1}{2})\}$.
Recall that
\begin{equation}
W^A=\Cn v_A,\  M^A=\C[e_0,\ld,e_{-n+1}]w_A.\notag
\end{equation}
Thus, to obtain $M^A$ from $W^A$ one must move $v_A$ to
$w_A$ and then "turn" $W^A$, to get $e_{1-i}$ from $e_{n-i}$.
One can check that while moving $v_A\rightsquigarrow w_A$ the
$z$-degree will decrease on $\sum \al_i$, and the $q$-degree will increase
on
$\sum \al_i^2-\frac{p(D)}{4}$ (note that in the formula (\ref{chW}) we
regard the q-degree of $v_A$ to be equal to zero, so
$\sum \al_i^2-\frac{p(D)}{4}$
is just a $q$-degree of $w_A$ in $F^{\otimes k}$).
Thus, to obtain from (\ref{chW}) the formula
for the character of $M^A$, one must multiply (\ref{chW}) on the
corresponding factor and replace $z$ by $\frac{z}{q^n}$, what is
corresponding to the "turn" of $W^A$
(we must divide on $q^n$, because the formula for the character of $W^A$
depends on the variables $(q,zq)$).
We obtain the following expression (substituting $n=\sum d_i$):
\begin{multline}
\ch(M^A,q,z)=
z^{-\Sum_{i=1}^k \al_i}
q^{\Sum_{i=1}^k \al_i^2} q^{-p(D)/4}\times\\
\times\sum_{j_k=0}^{d_{k+1}}\sum_{j_{k-1}=0}^{d_k+j_k}\ld
\sum_{j_1=0}^{d_2+j_2}
{\left(\frac{z}{q^{\sum_{i=2}^{k+1} d_i}}\right)}^{\sum\limits_{l=1}^k j_l}
q^{\sum\limits_{l=1}^k j_l(d_2+\cdots +d_l+j_l)}\times\\
\times\qb{d_{k+1}}{j_k}\qb{d_k+j_k}{j_{k-1}}\ld\qb{d_2+j_2}{j_1}.
\end{multline}
After a simple rearranging, we obtain:
\begin{multline}
\ch(M^A,q,z)=q^{-p(D)/4}\times\\
\times\sum_{j_k=0}^{d_{k+1}}\sum_{j_{k-1}=0}^{d_k+j_k}\ld
\sum_{j_1=0}^{d_2+j_2}
      z^{\Sum_{l=1}^k (j_l-\al_l)}
q^{\sum\limits_{l=1}^k (j_l-\al_l)^2}\times\\
\times\qb{d_{k+1}}{j_k}\qb{d_k+j_k}{j_{k-1}}\ld\qb{d_2+j_2}{j_1}.
\end{multline}
Make a change $i_l=j_l-\al_l$:
\begin{multline}
\ch(M^A,q,z)=q^{-p(D)/4}\times\\
\times\sum_{\gf
{i_k=-\al_k}{i_k+\al_k\in\Z}}^{d_{k+1}/2}
\sum_{\gf{i_{k-1}=-\al_{k-1}}{i_{k-1}+\al_{k-1}\in\Z}}^
                                    {d_k/2+i_k}\ld
\sum_{\gf{i_1=-\al_1}{i_1+\al_1\in\Z}}^{d_2/2+i_2}
      z^{\Sum_{l=1}^k i_l}
q^{\Sum_{l=1}^k i_l^2}\times\\
\times\qb{d_{k+1}/2+\al_k}{i_k+\al_k}
\qb{d_k/2+i_k+\al_{k-1}}{i_{k-1}+\al_{k-1}}
\ld\qb{d_2/2+i_2+\al_1}{i_1+\al_1}.
\end{multline}
Rewrite the binomial coefficients:
\begin{multline}
\label{chM}
\ch(M^A,q,z)=q^{-p(D)/4}
\sum_{\gf{i_l\ge -\al_l,\ i_l+\al_l\in\Z}
{i_k\le \frac{d_{k+1}}{2},\ i_{l+1}\ge i_l-\frac{d_{l+1}}{2}}}
   z^{\Sum_{l=1}^k i_l}
q^{\sum\limits_{l=1}^k i_l^2}\times\\
\times
\frac{1}{\qf{i_k+\frac{d_{k+1}}{2}}
\Prod_{l=1}^{k-1}\qf{i_{l+1}-i_l+\frac{d_{l+1}}{2}}}
\times\\
\times\frac{\qf{d_{k+1}}}{\qf{\frac{d_{k+1}}{2}+i_k}}
\Prod_{l=1}^{k-1} \frac{\qf{i_{l+1}+\frac{d_{l+1}}{2}+\al_l}}{\qf{i_l+\al_l}}
\end{multline}
We are interesting in the limit of (\ref{chM}) while $d_{k+1}\to\infty$
($d_{k+1}$ must preserve its parity). Denote
$$(\infty)_q=\Prod_{i=1}^{\infty} (1-q^i).$$
(In our formulas the symbol $(\infty)_q$ for the
$\prod_{i=1}^{\infty} (1-q^i)$
seems to be more natural, than the more common one $(q)_{\infty}$).
In order to compute the limit, note that we can consider such $i_l$ that
$-d_{k+1}/4<i_l<d_{k+1}/4$ (otherwise the $q$-degree will overreach
$d_{k+1}^2/16$).
Then the degree in $q$ of the polynomial
$$
\frac{\qf{i_{l+1}+\frac{d_{l+1}}{2}+\al_l}}{\qf{i_l+\al_l}}-1
$$
is greater than $\frac{d_{k+1}}{4}$.
So
$$
\frac{\qf{d_{k+1}}}{\qf{\frac{d_{k+1}}{2}+i_k}}
\Prod_{l=1}^{k-1}
\frac{\qf{i_{l+1}+\frac{d_{l+1}}{2}+\al_l}}{\qf{i_l+\al_l}}
\longrightarrow 1 \text{ while } d_{k+1}\to\infty.
$$
In the same time
$$
\frac{1}{\qf{i_k+\frac{d_{k+1}}{2}}}\longrightarrow \frac{1}{(\infty)_q}.
$$
We have proved the following theorem:
\begin{theorem}
\label{chL}
Let $D\in(\N\cup 0)^{k+1},\ d_{k+1}\le 1$. Denote
$2\al_i=d_{k+1}+\ld+d_{i+1}$,\ $i=1,\ld,k$. Then
$$
\ch L^D=\frac{q^{-p(D)/4}}{(\infty)_q}
\sum_{\gf{\gf{i_1,\ld,i_k}{i_l+\al_l\in\Z, l=1,\ld,k}}
     {i_l\ge i_{l-1}-\frac{d_l}{2}, l=2,\ld,k}}
\frac{ z^{\Sum_{l=1}^k i_l}q^{\sum\limits_{l=1}^k i_l^2}}
{\Prod_{l=1}^{k-1}\qf{i_{l+1}-i_l+\frac{d_{l+1}}{2}}}.
$$
\end{theorem}
\begin{cor} We obtain the formula from \cite{sto}:
\label{irrch}
$$
\ch L_{j,k}=\frac{1}{(\infty)_q}
\sum_{\gf{i_1,\ld,i_k\in\Z}
     {i_1\le\ld\le i_k}}
\frac{z^{\frac{j}{2}+\Sum_{l=1}^k i_l}
q^{\Sum_{l=1}^k i_l^2+\Sum_{l=1}^j i_l}}
{\Prod_{l=1}^{k-1} \qf{i_{l+1}-i_l}}.
$$
\end{cor}

\subsection{Combinatorial recurrent formula for $L^D$.}
Introduce the following notations:
\begin{gather}
D=(\ld d_{s-1} (d_s+2) d_{s+1}\ld), D'=(\ld d_{s-1} d_s d_{s+1}\ld),\\
D''=(\ld (d_{s-1}+1) d_s (d_{s+1}+1)\ld).\notag
\end{gather}
\begin{prop}
\label{decch}
$\ch L^{D}=\ch L^{D'}+q^{\frac{p(D'')-p(D)}{4}+\frac{d_s+1}{2}}
\ch L^{D''}.$
\end{prop}
\begin{proof}
Note that $p(D)=p(D').$ We have:
\begin{multline}
\label{form}
\ch L^{D}-\ch L^{D'}=\frac{q^{-p(D)/4}}{(\infty)_q}\times\\
\times\biggl(
\sum_{\gf{\gf{i_l+\al_l\in\Z, l=1,\ld,k}{i_s\ge i_{s-1}-\frac{d_s+2}{2}}}
  {i_l\ge i_{l-1}-\frac{d_l}{2}, l\ne s}}
\frac{z^{\Sum_{l=1}^k i_l}q^{\sum\limits_{l=1}^k i_l^2}}
{\qf{i_s-i_{s-1}+\frac{d_s+2}{2}}
\Prod_{l=1,l\ne s-1}^{k-1}\qf{i_{l+1}-i_l+\frac{d_{l+1}}{2}}}-\\
-\sum_{\gf{i_l+\al_l\in\Z, l=1,\ld,k}
     {i_l\ge i_{l-1}-\frac{d_l}{2}, l=2,\ld,k}}
z^{\Sum_{l=1}^k i_l}q^{\sum\limits_{l=1}^k i_l^2}
\frac{1}{\Prod_{l=1}^{k-1}\qf{i_{l+1}-i_l+\frac{d_{l+1}}{2}}}\biggr).
\end{multline}
Divide the summation region of the minuend into the two parts: the first one
with a parameters
with $i_{s+1}=i_s-\frac{d_s+2}{2}$, and the second one with
$i_{s+1}\ge i_s-\frac{d_s}{2}$.
Then the conditions on the parameters in the second part will coincide with
the conditions on the parameters in subtrahend. Thus the previous expression
equals
\begin{multline}
\frac{q^{-p(D)/4}}{(\infty)_q} \biggl(
\sum_{\gf{\gf{i_l+\al_l\in\Z, l=1,\ld,k}{i_s=i_{s-1}-\frac{d_s+2}{2}}}
  {i_l\ge i_{l-1}-\frac{d_l}{2}, l\ne s}}
 z^{\Sum_{l=1}^k i_l}q^{\sum\limits_{l=1}^k i_l^2}
\frac{1}{\Prod_{l=1,l\ne s-1}^{k-1}\qf{i_{l+1}-i_l+\frac{d_{l+1}}{2}}}+\\
+\sum_{\gf{i_l+\al_l\in\Z, l=1,\ld,k}
          {i_l\ge i_{l-1}-\frac{d_l}{2}, l=2,\ld,k}}
z^{\Sum_{l=1}^k i_l}q^{\sum\limits_{l=1}^k i_l^2}
\frac{q^{i_s-i_{s-1}+\frac{d_s+2}{2}}}
{\qf{i_s-i_{s-1}+\frac{d_s+2}{2}}
\Prod_{l=1,l\ne s-1}^{k-1}\qf{i_{l+1}-i_l+\frac{d_{l+1}}{2}}}\biggr).
\end{multline}
Add the first summand to the second:
\begin{equation}
\frac{q^{-p(D)/4}}{(\infty)_q}
\sum_{\gf{\gf{i_l+\al_l\in\Z, l=1,\ld,k}{i_s\ge i_{s-1}-\frac{d_s+2}{2}}}
  {i_l\ge i_{l-1}-\frac{d_l}{2}, l\ne s}}
\frac{z^{\Sum_{l=1}^k i_l}q^{\sum\limits_{l=1}^k i_l^2}
q^{i_s-i_{s-1}+\frac{d_s+2}{2}}}
{\qf{i_s-i_{s-1}+\frac{d_s+2}{2}}
\Prod_{l=1,l\ne s-1}^{k-1}\qf{i_{l+1}-i_l+\frac{d_{l+1}}{2}}}.
\end{equation}
Let us change the parameters $i_s:=i_s+1/2,\ i_{s-1}:=i_{s-1}-1/2$.
Then we obtain:
$$
\ch L^{D}-\ch L^{D'}=q^{\frac{p(D'')-p(D)+2d_s+2}{4}}\ch L^{D''}.
$$
Proposition is proved.
\end{proof}
\begin{cor}
$L^{D}=L^{D'}\oplus L^{D''}.$
\end{cor}

\newcounter{a}
\setcounter{a}{2}

\end{document}